\documentclass[12pt]{article}
\usepackage[english]{babel}
\usepackage{amsthm,amsmath,amssymb,amsfonts,latexsym}

\usepackage[usenames]{color}

\newtheorem{theorem}{\qquad Theorem}

\newtheorem{lemma}{\qquad Lemma}
\newtheorem{definition}{\qquad Definition}

\begin{document}
\newpage
%\tableofcontents %Table of contents
\cleardoublepage %The first chapter should start on an odd page.

\pagestyle{plain} %Now display headings: headings / fancy / ...

\title{On some polynomial version on the sum-product problem for subgroups}
\maketitle
\begin{center}
\textbf{\abstractname}
\end{center}

We generalize two results about subgroups of multiplicative group of finite field of prime order. In particular, the lower bound on the cardinality of the set of values of polynomial $P(x,y)$ is obtained under the certain conditions, if variables $x$ and $y$ belong to a subgroup $G$ of the multiplicative group of the filed of residues. Also the paper contains a proof of the result that states that if a subgroup $G$ can be presented as a set of values of the polynomial $P(x,y)$, where $x\in A$, and $y\in B$ then the cardinalities of sets $A$ and $B$ are close (in order) to a square root of the cardinality of subgroup $G$.

\section{Intoduction}

Let $\mathbb{F}_p=\mathbb{Z}/p\mathbb{Z}$ be a finite field of a prime order $p$, and $\mathbb{F}_p^*$ is its multiplicative group. 
Consider the polynomial $P\in\mathbb{F}_p[x,y]$. Let us define the set
\begin{eqnarray}\label{PAB}
P(A,B)=\{ P(a,b)\mid a\in A,\,b\in B \},
\end{eqnarray}
where $A$ и $B$ are subsets of $\mathbb{F}_p$, which can be called as {\it polynomial sum of sets $A$ and $B$}.  The particular case of such polynomial sum is the sum of sets
$$
A+B=\{ a+b \mid a\in A,\,b\in B\}.
$$
Let $G$ be a subgroup the group $\mathbb{F}_p^*$. In this part we consider the case $A=B=G$. For the cardinality of $|G+G|$ the following bounds have been obtained. As a corollary of the bound of~\cite{GV} for subgroup $G$ such as $|G|\ll p^{3/4}$ the following bound was obtained:
$$
|G\pm G|\gg |G|^{4/3}.
$$
In this formula and further symbols ``$\ll$'' и ``$\gg$'' are Vinogradov's symbols.

Heath-Brown and Konyagin proved the inequality (see~\cite{HK}) :
$$
|G\pm G|\gg |G|^{3/2}
$$
for subgroups $|G|\ll p^{2/3}$. 
The bound
$$
|G\pm G|\gg \frac{|G|^{5/3}}{\log^{1/2}|G|}.
$$
for such subgroups that $|G|\ll p^{1/2}$ is obtained in~\cite{VS}.

\bigskip

The second problem touches the possibility of presenting $G$ as a 
$$
G=P(A,B),
$$
where $P(A,B)$ is defined in (\ref{PAB}). Let $A$ and $B$ are non-trivial (sizes of $A$ and $B$ are exceed one) subsets of the set of residues modulo prime number $p$.
In the second part of the paper it is proved that if it is possible, then the cardinality of $|A|$ and $|B|$ are close to $\sqrt{|G|}$ (see part~\ref{pol-sumset}). This result generalizes the result of Shparlinski (see Th. 8 in~\cite{Shp}) to the polynomials $P(x,y)$ that are more general than $P(x,y)=x+y$.

\section{Polynomials on subgroups}\label{pol-sumprod}

\begin{definition} Let us call the polynomial $P\in \mathbb{F}_p[x,y]$ {\it good} if it is homogeneous with respect to $x$ and $y$, polynomial $P(x,y)-1$ is absolutely irreducible (it is irreducible over the algebraic closure $\overline{\mathbb{F}}_p$ of the field $\mathbb{F}_p$) and at least one of the polynomials $P(x,0)$, $P(0,y)$ is not identity to zero.
\end{definition}

\begin{definition}\label{def-subgr}
For a prime number $p$ and a natural number $n$ let us call a subgroup $G\subset \mathbb{F}_p^*$ $(n,p)$-admitted if
$$
100n^3<|G|<\frac{1}{3}p^{1/2}.
$$
\end{definition}

Theorem $2$ of the paper~\cite{VM} for homogeneous polynomial $P(x,y)$ can be re-formulated as follows.

\begin{theorem}\label{VM-th} 
For any $n$ there exist constants $C_1, C_2 > 0$ such that: for any prime $p$, $(n, p)$-admitted subgroup $G\in\mathbb{F}_p^*$, a good polynomial $P(x,y)$ of degree $n$, a natural number $h<C_2|G|^{2}$ and numbers $\alpha_1,\ldots,\alpha_h\in \mathbb{F}_p^*$ belonging to different $G$-cosets, there are at most
$$
C_1h^{2/3}|G|^{2/3}
$$
pairs $(x,y)$, for which $P(x,y)=\alpha_k$ for at least one $k=1,\ldots,h$. 
\end{theorem}

Values of constants 
$$
C_1=24n^4,\qquad C_2=40^{-3}n^{-9}
$$
were set in~\cite{VM}.
%Polynomials $P(x,y)-\alpha_k$ from Theorem \ref{VM-th} are also absolutely irreducible. 
Let us prove it in the following lemma.

\begin{lemma} 
If $P(x,y)$ is a good polynomial then the polynomial $P(x,y)-\alpha$, where $\alpha\in\mathbb{F}_p^*$ is also absolutely irreducible.
\end{lemma}

\textbf{Proof.} For any $\alpha\in\mathbb{F}_p^*$ let us denote by $a$ an arbitrary root of the $n$-th power of $1/\alpha$ in the algebraic closure $\overline{\mathbb{F}}_p$ ($a=\sqrt[n]{1/\alpha}$). Introduce the polynomial
$$
P_a(x,y)=P(ax,ay)-1,
$$
and suppose that the polynomial $P_a(x,y)$ is reducible
\begin{eqnarray}\label{P_a}
P_a(x,y)=P(ax,ay)-1=P_1(x,y)P_2(x,y).
\end{eqnarray} 
Let us substitute $x/a$ and $y/a$ instead of $x$ and $y$ into the equation (\ref{P_a}), then we obtain that
$$
P_a\left(\frac{x}{a},\frac{y}{a}\right)=P(x,y)-1=P_1\left(\frac{x}{a},\frac{y}{a}\right)P_2\left(\frac{x}{a},\frac{y}{a}\right),
$$ 
i.e. $P(x,y)-1$ is reducible. That contradicts to the assumption. So, we have that
$$
P_a(x,y)=P(ax,ay)-1=a^nP(x,y)-1=\frac{P(x,y)}{\alpha}-1
$$ 
is irreducible. Multiplying $P_a(x,y)$ by $\alpha$ there would be irreducible polynomial $P(x,y)-\alpha=\alpha P_a(x,y)$. $\square$

\begin{theorem}\label{nontriv-est}
For any $n$ there exists $C>0$ such that for any prime number $p$, $(n,p)$-admitted subgroup $G\in\mathbb{F}_p^*$ and a good polynomial $P(x,y)$ of degree $n$ we have the bound 
$$
|P(G,G)|>C|G|^{3/2}.
$$
\end{theorem}

\textbf{Proof.} Suppose the contrary. Then there exists such $n$ that  the statement of the theorem is not satisfied. That means that for any constant $C$, there are subgroup $G$ and a polynomial $P(x,y)$, with the given properties such that
$$
|P(G,G)|\le C|G|^{3/2}.
$$ 
Such pairs $(P,G)$ for the constant $C$ we call bad.

We apply Theorem~\ref{VM-th} to obtain the contradiction: for given $n$, it needs to be chosen $C_1,\, C_2 > 0$, satisfying the conditions of Theorem~\ref{nontriv-est}. After that let us put $C>0$ such that 
$$
C<C_2;\text{  } C_1C^{2/3}<\frac{100n^2-1}{100n^2}.
$$
The reason for put it like that will be clear later on. 

Let us take any bad pair $(P,G)$ for the chosen $C$. All possible values of $P(G, G)$ that are not greater than $C|G|^{3/2}$ and non-zero, can be arranged in the form of the Young tableau in such a way that each row contains values from one $G$-coset, and in different rows there are from different cosets. Thus, each line of the resulting diagram has no more than $|G|$ elements. Let us estimate from above the number of pairs $(x,y)$, for which the value lies into one or another column.

1) The number of pairs for which $P(x,y)=0$ is not greater than $n|G|$. 

Indeed, the polynomial $P(x,y)$ is homogeneous, that means that when $x=x_0\not= 0$ the polynomial $P(x_0,y)\in\mathbb{F}_p[y]$ is non-identity to zero. It has no more than $n$ roots. Let us estimate the number of pairs $(x,y)$ such that 
\begin{eqnarray}\label{pairsGG}
P(x,y)=0,\qquad (x,y)\in G\times G.
\end{eqnarray}
Let $x_0\in G$, that means that $x_0\not=0$, then the number of pairs $(x_0,y)\in G\times G$, $P(x_0,y)=0$ is not greater than $n$, therefore, the overall number of pairs (\ref{pairsGG} is not more than $n|G|$, since for every $x\in G$ there exist no more than $n|G|$ pairs.

2) If any column has $h$ elements then it can be noted that 
$$
h\leqslant|P(G, G)| \leqslant C|G|^{3/2}<C_2|G|^{3/2},
$$ 
therefore, since all the elements of the column lie in different cosets, according to Theorem~\ref{VM-th}, there exist at most $C_1h^{2/3}|G|^{2/3}$ pairs $(x,y)$ for which $P(x,y)$ lies into this column.

Now it can be denoted the column lengths for $h_1, h_2, \cdots, h_{|G|}$ and estimate the total number of pairs:
$$
|G|^2 < n|G|+\sum_{k=1}^{|G|}C_1h_k^{2/3}|G|^{2/3}.
$$

On the other hand, by the inequality on the power averages:
$$
\left(\frac 1{|G|}\sum_{k=1}^{|G|}h_k^{2/3}\right)^{3/2}\le \frac 1{|G|}\sum_{k=1}^{|G|}h_k.
$$
The sum of all $h_k$ is the total number of cells in the table, чso it does not exceed $C|G|^{3/2}$, whence:
$$
|G|^2<n|G|+C_1|G|^{2/3}\cdot|G|\left( \frac{C|G|^{3/2}}{|G|} \right)^{2/3} = n|G|+C_1C^{2/3}|G|^2 < n|G|+\frac{(100n^2-1)|G|^2}{100n^2}.
$$
As $|G|>100n^3$ (see Definition~\ref{def-subgr}), it is a contradiction, therefore, the theorem is proved. 

The value of the constant $C$ is following
$$
C=\min\left(\left(\frac{100n^2-1}{100n^2C_1}\right)^{3/2};\, C_2\right). %= \left(\frac{100n^2-1}{3200n^7}\right)^{3/2},
$$
%as $C_2$ was not found.
$\square$

\subsection{On additive shifts of multiplicative subgroups}

Using some algebraic ideas Garcia and Voloch in 1988 (see~\cite{GV})  proved that for any multiplicative subgroup $G\subseteq\mathbb{F}_p^{*}$ such that $|G|<(p-1)/((p-1)^{1/4}+1)$ and any non-zero $\mu$:
\begin{equation}
\label{trivial12}
|G\cap(G+\mu)|\leq 4|G|^{2/3}.
\end{equation}
Heath-Brown and Konyagin using the Stepanov's method (see~\cite{Step}) simplified the proof of this result and improved the constants in 2000 (see~\cite{HK}). In 2012 Vyugin and Shkredov generalize this bound to  the case of several additive shifts (see \cite{VS}).

Heath-Brown and Konyagin proved the inequality:
\begin{equation}
\label{trivial14}
|G\pm G|\gg |G|^{4/3}
\end{equation} 
for all subgroups $G$ for which $|G|\ll p^{2/3}$. 
Vyugin and Shkredov improved the inequality (\ref{trivial14}):
$$
|G\pm G|\gg \frac{|G|^{5/3}}{\log^{1/2}|G|}
$$
for subgroups $G$ such that $|G| \ll p^{1/2}$ (see~\cite{VS}).

\subsection{Polynomial version of sum-set problem}\label{pol-sumset}

Consider a subgroup $G\subset \mathbb{F}_p^* $, $G$-cosets $G_1,...,G_n$ ($G_i=g_iG$, where $g_i\in\mathbb{F}_p^*$, $1\le i\le n$ are arbitrary, they can be the same) and also consider the mapping 
$$
f: x \longmapsto (f_1(x),...,f_n(x))\in \mathbb{F}_p^n,\qquad n\ge 2 
$$
with polynomials $f_1(x),...,f_n(x)\in \mathbb{F}_p[x]$.

%\textbf{Definition 3.} {\it 
\begin{definition}
Let us call the set of polynomials $f_1(x),...,f_n(x)$ permissible if every polynomial $f_i(x)$ has at least one root $x_i\neq 0$ (in algebraic closure $\overline{\mathbb{F}}_p$ of the filed $\mathbb{F}_p$), which is not congruent with any of other roots of the polynomial set, that means 
$$
f_i(x_i)=0,\quad f_j(x_i)\neq 0,\,\, i \neq j, \,\,1 \le i,j\le n;\qquad x_i\not= x_j,\quad i\not= j
$$
and has non-zero free member $f_i(0)\neq 0$, $i=1,...,n$.
\end{definition}
%} 
In the paper \cite{V} there was obtained the higher estimation of the cardinality of set $M$:
$$
M=\{x\,|\, f_i(x)\in G_i, i=1,...,n\}.
$$

\begin{theorem}\label{th-map} Let $G$ be a subgroup of $\mathbb{F}_p^* $ ($p$ is a prime number), $G_1,...,G_n$ are $G$-cosets, $f_1(x),...,f_n(x)$ is a permissible set of polynomials of degrees $m_1,...,m_n$ respectively. Let the following inequality be true: 
$$
C_1(m,n) <|G|<C_2(m,n)p^{1-1/(2n+1)},
$$
where $C_1(m,n),\,C_2(m,n)$ are constants, depending on $n$ and $m=(m_1,...,m_n)$. Then the following estimate
$$
|M| \le C_3(m,n)|G|^{1/2+1/(2n)}
$$
is correct. Constants can be chosen as follows:
$$
C_1(m, n)=2^{2n}(\max m_i)^{4n},
\qquad
C_2(m, n)=(n+1)^{-\frac{2n}{2n+1}}(m_1\ldots m_n)^{-\frac{2}{2n+1}},
$$
$$
C_3(m, n)=4(n+1)(m_1\dots m_n)^{\frac{1}{n}}\sum_{i=1}^n m_i.
$$
\end{theorem}

\begin{definition}
Let us call the polynomial $P(x,y)\in \mathbb{F}_p[x,y]$ required if it cannot be divided by any of polynomials of neither $x$ or $y$, except constants, that means
$$
f(x)\mid P(x,y)\Rightarrow f(x)\equiv {\rm const};
$$
$$
g(y)\mid P(x,y)\Rightarrow g(y)\equiv {\rm const}.
$$
\end{definition}

\begin{lemma}\label{lem-dop} 
For any required polynomial $P(x,y)$, where $\deg_x P=k$, $\deg_y P=l$, among polynomials $f_i(x)=P(x,y_i)$, where $y_i$ are different elements of $\mathbb{F}_p$, $i=1,...,h$, there can be found the permissible subset $f_{i_1},...,f_{i_N}$ from $N=\left[\frac{h-2l}{kl}\right]$ polynomials.
\end{lemma}

\textbf{Proof.} It can be noted that the number $x=r$ can be the root at most of $l$ polynomials $f_i(x)=P(x, y_i)$. The contrary would mean that the polynomial $g(y)=P(r,y)$ has more than $l$ roots, but its degree is not greater than $\deg_y P(x,y)=l$. Therefore, it has to be zero but in this case $P(x, y)$ is divided by $(x-r)$, that contradicts the fact that $P(x,y)$ is required.

Firstly, let us take out from the set $y_1,\ldots, y_h$ all such $y_i$ that are the roots of the leading coefficient а $p_k(y)$ and a free term $p_0(y)$
of polynomial
$$
P(x,y)=p_k(y)x^k+\ldots+p_0(y),
$$ 
being considered as a polynomial of the variable $x$. It is obvious that the number of roots is not greater than $2l$, as both leading and free terms are non-zero polynomials of variable $y$, which degree is not greater than $l$ (free term is non-zero as $P$ is required and cannot be divided by $x$).

From remaining not less than $(h-2l)$ values $y_i$ it can be chosen any: $P(x, y_i)$ has no more than $k$ roots (as the leading term is non-zero). Let us take out all $y_j$ such that $P(x, y_j)$ has at least one common root with $P(x, y_i)$. From above it can be see that for every polynomial $P(x, y_i)$ that has no more than $k$ roots there exist no more than $l$ polynomials from the set, that have this as a root. Therefore, there are no more than $kl$ polynomials, that have common root with $P(x, y_i)$. Let us repeat this process: from remaining $y_i$ it can be chosen one and taken out no more than $kl$ values $y_j$ such that this polynomial has at least one common root with the considered polynomial. At the end it can be chosen minimum $\left[\frac{h-2l}{kl}\right]$ polynomials $P(x, y_i)$, none of two of each have no common roots. Also it can be seen that these polynomials have non-zero free term as it was taken out all $y_i$ that make it zero. Therefore the taken set is permissible. $\square$

\begin{theorem}\label{th-G12}%\textbf{Теорема 5.} {\it 
For any $k$ and $l$ there can be found the constant $C(k, l)$ such as for any $G\subset\mathbb{F}_p^*,$ required polynomial $P(x,y)$ of degrees $k$ and $l$ on $x$ and $y$ respectively, $A,B\subset \mathbb{F}_p$ with conditions 
$$
|G|<C p^{1-o(1)},
$$ 
$$
G=P(A,B),
$$ 
$$
|A|,|B|\gg 1,
$$  
the cardinalities of sets $A$ and $B$ are of order $|G|^{1/2+o(1)}$.
%}
\end{theorem}

\textbf{Proof.} Let $h(n, k, l)$  be the minimum $h$, which has to be taken in Lemma~\ref{lem-dop}, so that from the set of $h$ values of $y$ there would be $n$ permissible polynomials. It exists and no more than $nkl+2l$ in Lemma~\ref{lem-dop}. Let $\delta>0$ and $\varepsilon>0$ be the indexes, which can be taken in the statement of the theorem instead of $o(1)$. It means that
$$
|G|<C p^{1-\delta},
$$ 
and it has to be proved that 
$$
|G|^{1/2-\varepsilon}<|A|,|B|<|G|^{1/2+\varepsilon}.
$$
Let us take $q\ge 2$ such that 
$$
1-1/(2q+1) > 1-\delta,
$$ 
and choose $C$ such that for every $p$:  
$$
C p^{1-\delta}<(p/k)^{1-1/(2q+1)}/(q+1).
$$
Let $|A|,|B|>h(q, k, l)$. Then due to Lemma~\ref{lem-dop} from $|B|$ values of $y$ it can be chosen $q$ such that if it substitutes in $P$, there would be the permissible set of $q$ polynomials. Let us apply  Theorem~\ref{th-map} to this set and cosets $G_i=G$, $i=1,\ldots,h$. It can be done since the last inequality will be transformed to 
$$
|G|<(p/k)^{1-1/(2q+1)}/(q+1),
$$ 
that follows from the first condition and choice of $q,\, C$. The constant in Theorem~\ref{th-G12} depends only on $k$ and $\delta$ as $m=(\underbrace{k,...,k}_{\text{$q$ times}})$. Left inequalities in Theorem~\ref{th-G12} are satisfied if $G$ is sufficiently large and $k, l, \delta$ are fixed. For small $G$ there is nothing to prove as 
$$
|A|,|B|\gg 1,\quad G=P(A,B).
$$ 
The set $M$ for such small cosets includes $A$. That means that
$$
|A|\le C_1(k, \delta)|G|^{1/2+1/(2q)} \le C_1(k, \delta)|G|^{3/4}.
$$ 
Applying the fact that
$$
|A||B|\ge |G|,
$$ 
as use of polynomial $P$ is a surjective mapping $A\times B\to G,$ then
$$
|B|\ge (1/C_1(k, \delta ))|G|^{1/4}.
$$
Hence, it can be proved that for any $n$ there exists constant $C_2(k, l, n, \delta  )$ such that
$$
|A|<C_2(k, l, n, \delta )|G|^{1/2+1/(2n)}.
$$ 
If
$$
(1/C_1(k, \delta))|G|^{1/4}\ge h(n, k, l),
$$ 
then from 
$$
|B|>h(q, k, l)
$$ 
it follows 
$$
|B|>h(n, k, l).
$$ 
Applying Theorem~\ref{th-map} one more time, for set of $n$ substitutions $y$ from $B$ and cosets, that are equal to $G$, then 
$$
|A|\le C_2(k,l,n, \delta )|G|^{1/2+1/(2n)}
$$ 
for every 
$$
|G|\ge(h(n, k, l)/C_1(k, \delta ))^4.
$$ 
The right part of the last inequality depends only on $k, l, n, \delta $, so increasing constant even more  $C_2(k,l,n, \delta )$, it can be obtained in other cases as well.

The same time, it can be obtained
$$
|B|\le C_3(k,l,n, \delta )|G|^{1/2+1/(2n)},
$$ 
using the symmetric condition. From 
$$
|A||B|\ge |G|
$$ 
it follows that for another constant $C_4(k, l, n, \delta)$  
$$
|A|,|B|\ge C_4(k,l,n, \delta )|G|^{1/2-1/(2n)}.
$$ 
As $n$ can be big as much as possible, $1/(2n)$ can be taken less than $\varepsilon$. The existence of such constants means that
$$
|G|^{1/2-\varepsilon}<|A|,|B|<|G|^{1/2+\varepsilon}.\ \ \ \square
$$
The authors are grateful to Andrey Volgin for his useful comments.

Sofia Aleshina\\
National Research University "Higher school of economics",\\
University of Bedfordshire,\\
{\it aleshina.sofia@mail.ru}.
\bigskip\\
Ilya Vyugin\\
The Institute for Information Transmission Problems of Russian Academy of Sciences,\\
National Research University "Higher school of economics",\\
Steklov Mathematical Institute of Russian Academy of Sciences,\\
{\it vyugin@gmail.com}.\\


\begin{thebibliography}{99}



\bibitem{VM} {\sc S. Makarychev, I. Vyugin}, {\em Solutions of polynomial equation over $\mathbb{F}_p$ and new bounds of
additive energy} // Arnold Math J. (2019) V 5, Issue 1, 105-121. 
   


\bibitem{Shp} {\sc I. E. Shparlinski} {\em Additive Decompositions of Subgroups of Finite Fields}, SIAM
J. Discrete Math., 27:4 (2013), 1870-1879.


\bibitem{GV} {\sc A. Garcia, J.Voloch}, {\em Fermat curves over nite elds}, J. Number Theory,
30:3 (1988), 345--356.

\bibitem{KS} {\sc S. V. Konyagin, I. E. Shparlinski}, {\em Character sums with exponential functions and their
applications}, Cambridge Tracts in Math., 136, Cambridge Univ. Press, Cambridge, 1999.


\bibitem{HK} {\sc D. Heath-Brown, S. Konyagin}, {\em New bounds for Gauss sums derived from k-thpowers, and for Heilbronn’s exponential sum}, Q. J. Math., 51:2 (2000), 221–235.


\bibitem{Step} {\sc S. Stepanov}, {\em On the number of points of a hyperelliptic curve over a finite prime field}, Math. USSR-Izv., 3:5 (1969), 1103–1114.


\bibitem{VS}
%{\sc Вьюгин~И.~В., Шкредов~И.~Д.,}
%{\em Об аддитивных сдвигах мультипликативных подгрупп, }
%Матем. сб. 203:{\bf 6} (2012), 81--100.
{\sc I. Vyugin, I. Shkredov}, 
{\em On additive shifts of multiplicative subgroups,} 
Sb. Math., 203:6 (2012), 844–863 


\bibitem{V} 
%{\sc И. В. Вьюгин}, {\em Оценка числа прообразов полиномиального отображения} //  Матем. заметки, 106:2 (2019), %212-221.
{\sc V’yugin, I.V.}, {\em A Bound for the Number of Preimages of a Polynomial Mapping} // Math Notes, 106, 203-211 (2019). 

\bibitem{CZ}
{\sc P. Corvaja, U. Zannier}, {\em Greatest common divisor of $u-1$, $v-1$ in positive characteristic and rational points on curves over finite fields} // J. Eur. Math. Soc., 15, 1927--1942, (2013).
345--356.
	
%\end{enumerate}
\end{thebibliography}
\end{document}